\documentstyle[12pt]{article}
\newcommand{\oid}[3]{\mbox{${\cal #1}_{#2}^{#3}$}}
\newcommand{\idn}[3]{\mbox{${\bf #1}_{#2}^{#3}$}}
\newcommand{\idrel}[3]{\mbox{${\cal #1} \stackrel{1}{#2} {\cal #3}$}}

\begin{document}
\newtheorem{theorem}{Theorem}[section]
\newtheorem{prop}{Proposition}[section]
\newtheorem{lemma}{Lemma}[section]
\newtheorem{cor}{Corollary}[section]

\title{Local properties of accessible injective operator ideals}
\author
{Frank Oertel\\
Zurich}
\maketitle

\begin{abstract}
{\noindent}In addition to Pisier's counterexample of a non-accessible maximal Banach ideal,  
we will give a large class of maximal Banach ideals which {\it{are accessible}}. 
The first step is implied by the observation that a "good behaviour" of
trace duality, which is canonically induced by conjugate operator ideals
can be extended to adjoint Banach ideals, if and only if these adjoint 
ideals satisfy an accessibility condition (theorem 3.1). This observation
leads in a natural way to a characterization of accessible injective
Banach ideals, where we also recognize the appearance of the ideal of 
{\it{absolutely summing operators}} (prop. 4.1). 
By the famous {\it{Grothendieck inequality}}, every operator from $L_1$ to 
a Hilbert space is absolutely summing, and therefore our search for such
ideals will be directed towards Hilbert space factorization - via an 
operator version of Grothendieck's inequality (lemma 4.2). As a consequence,
we obtain a class of injective ideals, which are "quasi-accessible", and with
the help of {\it{tensor stability}}, we improve the corresponding norm
inequalities, to get accessibility (theorem 4.1 and 4.2). In the last chapter
of this paper we give applications, which are implied by a non-trivial 
link of the above mentioned considerations to normed products of operator ideals. \\     

{\noindent}{\it{Key words and phrases:}} accessibility, Banach spaces, conjugate 
operator ideals, Hilbert space factorization, Grothendieck's inequality,  
tensor norms, tensor stability\\

{\noindent}{\it{1991 AMS Mathematics Subject Classification:}} primary 46M05, 47D50; 
secondary 47A80.  
\end{abstract}

\section{Introduction}

Given Banach spaces $E$, $F$ and a maximal Banach ideal (\oid{A}{}{}, \idn{A}{}{}),
we are interested in reasonable sufficient conditions on $E$, $F$ and
(\oid{A}{}{}, \idn{A}{}{}) such that (\oid{A}{}{}, \idn{A}{}{})
is accessible. In general it is a nontrivial subject to prove accessibility of  
maximal Banach ideals since non-accessibility can only appear on Banach
spaces without the metric approximation property, and in 1992, Pisier
made use of such a Banach space (the Pisier space $P$) to construct a non-accessible
maximal Banach ideal (cf. \cite{df}, 31.6). On the other hand, accessible Banach ideals
allow a suggestive (algebraic) calculus which leads to further results concerning
the local structure of operator ideals (e.g. a transfer of the principle of
local reflexivity from the operator norm to suitable
ideal norms \idn{A}{}{} (cf. \cite{df},\cite{oe1} and \cite{oe2}).\\

{\noindent}This paper is mainly devoted to the description of a large class of maximal injective Banach ideals
which are totally accessible. We will see a deep interplay between conjugates of 
Banach ideals, Hilbert space factorization, Grothendieck's inequality and 
tensor stable quasi-Banach ideals.\\
We only deal with Banach spaces and most of our notations and definitions 
concerning Banach spaces and operator ideals are standard and can be
found in the detailed monographs \cite{df} and \cite{p1}. However, if 
$(\oid{A}{}{}, \idn{A}{}{})$ and $(\oid{B}{}{}, \idn{B}{}{})$ are given
quasi-Banach ideals,
we will use the shorter notation $(\oid{A}{}{d}, \idn{A}{}{d})$ for the dual
ideal (instead of $(\oid{A}{}{dual}, \idn{A}{}{dual})$) and the abbreviation 
\idrel{A}{=}{B} for the equality 
$(\oid{A}{}{}, \idn{A}{}{}) = (\oid{B}{}{}, \idn{B}{}{})$. The inclusion 
$(\oid{A}{}{}, \idn{A}{}{}) \subseteq (\oid{B}{}{}, \idn{B}{}{})$ is often
shortened 
by \idrel{A}{\subseteq}{B}, and if $T : E \longrightarrow F$
is an operator, we indicate that it is a metric injection by writing 
$T : E \stackrel{1}{\hookrightarrow} F$.\\
{\noindent}Each section of this paper includes the more special terminology which
is not so common.

\section{On tensor norms and associated Banach ideals}
At first we recall the basic notions of Grothendieck's metric theory of tensor products
(cf., eg., \cite{df}, \cite{gl}, \cite{gr}, \cite{l}), which will be used throughout this
paper.\\
A {\it{tensor norm}} $\alpha$ is a mapping which assigns to each pair $(E, F)$
of Banach spaces a norm $\alpha(\cdot; E, F)$ on the algebraic tensor product
$E \otimes F$ (shorthand: $E {\otimes}_\alpha F$  and $E \tilde{\otimes}_\alpha F$
for the completion) such that
\begin{enumerate}
\item[(i)]   $\varepsilon\leq\alpha\leq\pi$
\item[(ii)]  $\alpha$ satisfies the metric mapping property: If
             $S \in {\cal L}(E, G)$ and $T \in {\cal L}(F, H)$,\\
             then $\| S\otimes T:E\otimes_\alpha F \longrightarrow G\otimes_\alpha H\|                    
             \leq \| S\|\ \| T\|$             
\end{enumerate}                          
Well-known examples are the injective tensor norm $\varepsilon$, which is the smallest one,
and the projective tensor norm $\pi$, which is the largest one. For other important
examples we refer to \cite{df}, \cite{gl}, or \cite{l}.\\
Each tensor norm $\alpha$ can be extended in two natural ways. For this, denote
for given Banach spaces $E$ and $F$
\begin{center}
FIN$(E) : = \{M \subseteq E\mid M\in $ FIN$\}$ \hspace{0.2cm}and\hspace{0.2cm}
COFIN$(E) : = \{L\subseteq E\mid E/L \in $ FIN$\}$,       
\end{center}
where FIN stands for the class of all finite dimensional Banach spaces.\\
Let $z \in E \otimes F$. Then the {\it{finite hull}}\/ $\stackrel{\rightarrow}{\alpha}$
is given by
\begin{center}
\hspace{0.2cm}$\stackrel{\rightarrow}{\alpha}(z; E, F) : = \inf\{\alpha(z; M, N) \mid M \in $ FIN$(E), 
N \in $ FIN$(F), z \in M \otimes N\}$
\end{center}
and the {\it{cofinite hull}}\/ $\stackrel{\leftarrow}{\alpha}$ of $\alpha$ is
given by
\begin{center}
\hspace{0.2cm}$\stackrel{\leftarrow}{\alpha}(z; E, F) : = \sup\{\alpha(Q_K^E 
\otimes Q_L^F(z); E/K, F/L) \mid K \in $ COFIN$(E), 
L \in $ COFIN$(F)\}$.
\end{center}
$\alpha$ is called {\it{finitely generated}}\/ if $\alpha = \hspace{0.1cm}
\stackrel{\rightarrow}{\alpha}$,
{\it{cofinitely generated}}\/ if $\alpha = \hspace{0.1cm}\stackrel{\leftarrow}{\alpha}$ 
(it is always true that $\stackrel{\leftarrow}{\alpha} \hspace{0.1cm}\leq \alpha \leq \hspace{0.1cm}\stackrel{\rightarrow}{\alpha}$).
$\alpha$ is called {\it{right-accessible}} if  
$\stackrel{\leftarrow}\alpha$$(z; M, F) = \hspace{0.1cm}
\stackrel{\rightarrow}\alpha$$(z; M, F)$ for all $(M, F) \in$
FIN $\times$ BAN, {\it{left-accessible}} if 
$\stackrel{\leftarrow}{\alpha}$$(z; E, N) = \hspace{0.1cm}
\stackrel{\rightarrow}{\alpha}$$(z; E, N)$ for all $(E, N) \in$
BAN $\times$ FIN, and {\it{accessible}} if it is right- and left-accessible.
$\alpha$ is called {\it{totally accessible}} if $\stackrel{\leftarrow}{\alpha}
\hspace{0.1cm}= \hspace{0.1cm}\stackrel{\rightarrow}{\alpha}$.\\ 
The injective norm $\varepsilon$ is totally
accessible, the projective norm $\pi$ is accessible - but not totally accessible,
and Pisier's counterexample implies the existence of a (finitely generated)
tensor norm which is neither left- nor right-accessible (see \cite{df}, 31.6).\\
There exists a powerful one-to-one correspondence between finitely generated
tensor norms and maximal Banach ideals which links thinking in terms of 
operators with "tensorial" thinking and which allows to transfer notions in
the "tensor-language" to the "operator-language" and conversely. We refer the 
reader to \cite{df} and \cite{oe1} for detailed informations
concerning this subject.\\
Let $E, F$ be Banach spaces and $z = \sum\limits_{i = 1}^n a_i \otimes y_i$\ be
an Element in $E' \otimes F$. Then $T_z(x): = \sum\limits_{i = 1}^n 
\langle x, a_i\rangle  y_i $
defines a finite operator $T_z \in \oid{F}{}{}(E, F)$ which is independent of the
representation of $z$ in $E' \otimes F$. 
Let $\alpha$ be a finitely generated tensor norm and (\oid{A}{}{}, \idn{A}{}{})
be a maximal Banach ideal. $\alpha$ and (\oid{A}{}{}, \idn{A}{}{}) are said to
be {\it{associated}}, notation:
\begin{center}
$(\oid{A}{}{}, \idn{A}{}{}) \sim \alpha$ \hspace{0.15cm}(shorthand: $\oid{A}{}{} 
\sim \alpha$, resp. 
\hspace{0.1cm}$\alpha \sim \oid{A}{}{}$)
\end{center}
if for all $M, N \in $ FIN 
\begin{center}
$\oid{A}{}{}(M, N) = M'{\otimes}_\alpha N$
\end{center}
holds isometrically: $\idn{A}{}{}(T_z) = \alpha(z; M', N)$.\\

{\noindent}Besides the maximal Banach ideal $(\cal L, \| \cdot \|) \sim \varepsilon$
we will mainly be concerned with $(\cal I, \bf I) \sim \pi$ (integral operators),
$(\oid{L}{2}{}, \idn{L}{2}{}) \sim w_2 $ (Hilbertian operators),
$(\oid{D}{2}{}, \idn{D}{2}{}) \stackrel{1}{=} 
(\oid{L}{2}{\displaystyle\ast}, \idn{L}{2}{\displaystyle\ast}) \sim w_2^\ast $\/
($2$-dominated operators),
$(\oid{P}{p}{}, \idn{P}{p}{}) \sim g_p\backslash = g_q^\ast $ (absolutely 
$p$-summing operators),
$1 \leq p \leq \infty, \frac{1}{p} + \frac{1}{q} = 1$,
$(\oid{L}{\infty}{}, \idn{L}{\infty}{}) \stackrel{1}{=} 
(\oid{P}{1}{\displaystyle\ast}, \idn{P}{1}{\displaystyle\ast}) \sim w_\infty $
and $(\oid{L}{1}{}, \idn{L}{1}{}) \stackrel{1}{=} 
(\oid{P}{1}{\displaystyle\ast d}, \idn{P}{1}{\displaystyle\ast d}) \sim w_1 $.\\
Since it is important for us, we recall the notion 
of the conjugate operator ideal (cf. \cite{glr}, \cite{jo}): \\
let $(\oid{A}{}{}, \idn{A}{}{})$
be a quasi-Banach ideal. Let $\oid{A}{}{\Delta}(E, F)$ be the set of all $T \in 
{\cal L}(E, F)$
for which 
\begin{center}
$\idn{A}{}{\Delta}(T) : = \sup\{tr( TL) \mid L \in {\cal F}(F, E), \idn{A}{}{}(L) \leq 1 \} 
< \infty $.
\end{center}
Then a Banach ideal is obtained. It is called the {\it{conjugate ideal}} of
$(\oid{A}{}{}, \idn{A}{}{})$.\\

{\noindent}$(\oid{A}{}{}, \idn{A}{}{})$ is called {\it{right-accessible}}, if for all
$(M, F) \in$ FIN $\times$ BAN, operators $T \in {\cal L}(M, F)$ and $ \varepsilon > 0$ 
there are $N \in$ FIN$(F)$ and $S \in {\cal L}(M, N)$ such that $ T = J_N^F S$
and $\idn{A}{}{}(S) \leq (1 + \varepsilon) \idn{A}{}{}(T) $.
It is called {\it{left-accessible}}, if for all $(E, N) \in$ BAN $\times$ FIN,
operators $T \in {\cal L}(E, N)$ and $\varepsilon > 0$ there are $L \in$ COFIN$(E)$
and $S \in {\cal L}(E/L, N)$ such that $T = S Q_L^E $ and 
$\idn{A}{}{}(S) \leq (1 + \varepsilon) \idn{A}{}{}(T) $.
A left- and right-accessible ideal is called {\it{accessible}}. 
$(\oid{A}{}{}, \idn{A}{}{})$ is {\it{totally accessible}}, if for every finite rank operator
$T \in {\cal F}(E, F)$ between Banach spaces and $\varepsilon > 0$ there are
$(L, N) \in$ COFIN$(E) \times$ FIN$(F)$ and $S \in {\cal L}(E/L, N)$ such that 
$T = J_N^F S Q_L^E$ and $\idn{A}{}{}(S) \leq (1 + \varepsilon) \idn{A}{}{}(T) $. \\
Every injective quasi-Banach ideal is right-accessible (every surjective ideal is
left-accessible) and, if it is left-accessible, it is totally accessible.\\
A finitely generated tensor norm is right-accessible (resp. left-accessible,
accessible, totally accessible) if and only if its associated maximal Banach ideal is.\\

\section{Accessible conjugate operator ideals}

Let (\oid{A}{}{}, \idn{A}{}{}) be a p-Banach ideal $(0 < p \leq 1)$.\footnote{
Most of the results in this paragraph first appeared in the author's 
doctoral thesis (see \cite{oe1}). However, we now are using different proofs
which give a better insight into the underlying structures.} 
Suppose \oid{A}{}{} is right-accessible. If we apply the cyclic composition theorem
(see \cite{df}, 25.4) to  \idrel{ A \circ L }{\subseteq}{A}, it follows that
$\oid{ A^{\displaystyle\ast} \circ A }{}{} \stackrel{1}{\subseteq} \oid{I}{}{}.$
If \oid{A}{}{} is totally accessible,
an easy calculation shows that \idrel{ A^{\displaystyle\ast}}{=}{A^{\scriptstyle\triangle}}.
For $p = 1$, these properties of \oid{A}{}{} characterize accessibility in the 
following sense:

\begin{theorem} Let $(\oid{A}{}{}, \idn{A}{}{})$ be a Banach
ideal. Then $(\oid{A}{}{{\displaystyle\ast}{\scriptstyle\triangle}}, \idn{A}{}{{\displaystyle\ast}{\scriptstyle\triangle}})$
is always right-accessible. $(\oid{A}{}{}, \idn{A}{}{})$ is right-accessible if and
only if \idrel{A^{\displaystyle\ast} \circ A}{\subseteq}{I}. If in addition 
$(\oid{A}{}{}, \idn{A}{}{})$ is maximal then $(\oid{A}{}{}, \idn{A}{}{})$ is 
totally accessible if and only if \idrel{ A^{\displaystyle\ast}}{=}{A^{\scriptstyle\triangle}}.
\end{theorem}

{\noindent}{\sc PROOF}:{\hspace{0.25cm}}To prove the right-accessibility of
\oid{A}{}{{\displaystyle\ast}{\scriptstyle\triangle}}, we may assume that
\oid{A}{}{} is maximal (cf \cite{p1}, 9.3). Let $\alpha \sim \oid{A}{}{}$ be
associated and $(M, F) \in$ FIN $\times$ BAN. Then $\alpha^{\displaystyle\ast} \sim \oid{A^{\displaystyle\ast}}{}{}$.
The representation theorem for minimal operator ideals (see \cite{df}, 22.2) gives

\begin{center}
$\oid{A}{}{min}(M, F)  \stackrel{1}{=} M' {\otimes}_\alpha F \stackrel{1}{\hookrightarrow}
(F {\otimes}_{{\alpha}^{t}} M')''$.
\end{center}
Since $\alpha$ is finitely generated, the representation theorem for maximal operator
ideals (see \cite{df}, 17.5) yields

\begin{center}
$(F {\otimes}_{{\alpha}^{t}} M')'\cong \oid{A}{}{\displaystyle\ast}(F, M)$.
\end{center}
On the other hand, by canonical trace duality, it follows that 

\begin{center}
$\oid{A}{}{\displaystyle\ast\scriptstyle\triangle}(M, F) \stackrel{1}{\hookrightarrow}
(\oid{A}{}{\displaystyle\ast}(F, M))'$.

\end{center}
Hence $\oid{A}{}{min}(M, F)  \stackrel{1}{=} 
\oid{A}{}{\displaystyle\ast\scriptstyle\triangle}(M, F)$.
Since \oid{A}{}{min} always is right-accessible (see \cite{df}, 25.3), 
\oid{A}{}{\displaystyle\ast\scriptstyle\triangle} is right-accessible.\\
Let (\oid{A}{}{}, \idn{A}{}{}) be an arbitrary Banach ideal such that 
\idrel{A^{\displaystyle\ast} \circ A}{\subseteq}{I}. First we show that for
all $(M, F) \in$ FIN $\times$ BAN, operators $T \in \oid{L}{}{}(M, F)$
\begin{center}
$\idn{A}{}{\displaystyle\ast\scriptstyle\triangle}(T) \leq \idn{A}{}{}(T)$
\end{center}
Let $L \in \oid{F}{}{}(F, M)$. Then $\mid tr(TL)\mid = \mid tr(LTid_M)\mid
\leq \idn{I}{}{}(LT) \cdot {\| id_M \|} \leq \idn{A^{\displaystyle\ast}}{}{}(L)
\cdot \idn{A}{}{}(T) \cdot 1$. Hence 
$\idn{A}{}{\displaystyle\ast\scriptstyle\triangle}(T) \leq \idn{A}{}{}(T)$. Since
\oid{A}{}{} is normed, $\idn{A}{}{}(S) = \idn{A}{}{\displaystyle\ast\displaystyle\ast}(S)
=  \idn{A}{}{\displaystyle\ast\scriptstyle\triangle}(S)$ for all {\it{elementary}}
operators $S$ - between finite dimensional spaces - (cf. \cite{p1}, 9.2.2), and it 
follows that \oid{A}{}{} is right-accessible.\\
Now let (\oid{A}{}{}, \idn{A}{}{}) be a Banach ideal such that it is maximal and 
\idrel{ A^{\displaystyle\ast}}{=}{A^{\scriptstyle\triangle}}. Let $\alpha \sim 
\oid{A}{}{}$ be associated, $E$, $F$ be Banach spaces and $z \in
E \otimes F$. Let $w: = j_E \otimes id_F(z)$ and $ T_w$ the associated operator
in $\oid{F}{}{}(E', F)$. Since $\alpha$ is finitely generated, the above
mentioned representation theorem for maximal operator ideals and a simple
application of the Hahn-Banach theorem give

\begin{eqnarray}
   \alpha(z; E, F) & = & {\alpha}^{t}(z^t; F, E) \nonumber\\
                   & = & \sup\{|\langle z^t,\phi\rangle| \mid \phi \in 
                          B_{(F \otimes_{\alpha^{t}} E)'}\} \nonumber\\
                   & = & \sup\{|tr(ST_w)| \mid S \in B_{\oid{A}{}{\displaystyle\ast}
                          \displaystyle(F, E')}\} \nonumber\\
                   & = & \sup\{|tr(ST_w)| \mid S \in B_{\oid{A}{}{\scriptstyle\triangle}
                          \displaystyle(F, E')}\} \nonumber
\end{eqnarray}
Hence  $\alpha(z; E, F) \leq \idn{A}{}{}(T_w) = \hspace{0.1cm}\stackrel{\leftarrow}
{\alpha}$$(z; E, F)$
(this equality follows from the embedding lemma (see \cite{df}, 17.6)).
Therefore $\alpha \sim \oid{A}{}{}$ is totally accessible, and the proof is
finished. ${}_{\displaystyle\Box}$\\

{\noindent} Given an arbitrary maximal Banach ideal (\oid{A}{}{}, \idn{A}{}{}) we 
have shown
that \oid{A}{}{\scriptstyle\triangle} is right-accessible. The natural 
question whether \oid{A}{}{\scriptstyle\triangle} is {\it{left}}-accessible
is still open and leads to interesting results concerning the local structure of
\oid{A}{}{}.\footnote{
Note that we cannot use the preceeding proof to verify the left-accessibility of
the conjugate ideal, since for all $(M, F) \in$ FIN $\times$ BAN, we have
$(M {\otimes}_{{\alpha}^{t}} F')'\cong \oid{A}{}{\displaystyle\ast}(M, F'')$.}\\ 
It is even true that \oid{A}{}{\scriptstyle\triangle} is left-accessible
{\it{if and only if}} the weak \oid{A}{}{}-local principle of 
reflexivity holds (i.e., in this case it is possible to transfer the 
estimation  for the operator norm $\| \cdot \|$  to 
the ideal norm \idn{A}{}{} (see \cite{oe1} and \cite{oe2} for further details)).\\
The dual ideal \oid{A}{}{d} is also a maximal Banach ideal and therefore
\oid{A}{}{d \scriptstyle\triangle} is right-accessible. Since 
$\oid{A}{}{\scriptstyle\triangle d} \stackrel{1}{\subseteq} 
\oid{A}{}{d \scriptstyle\triangle}$
we obtain that both \oid{A}{}{\scriptstyle\triangle d} and 
\oid{A}{}{\scriptstyle\triangle dd} 
are accessible. These considerations imply a slight generalization of
(\cite{df}, 25.11) which does not assume accessibility-conditions:
\begin{prop}
Let $(\oid{A}{}{}, \idn{A}{}{})$ be Banach ideal and $E, F$ be Banach spaces.\\
If $E'$ or $F$ has the approximation property, then
\begin{center}
$({\cal A}^{sur})^{min}(E, F) \stackrel{1}{=} ({\cal A}^{min})^{sur}(E, F)$
\end{center}
and
\begin{center}
$({\cal A}^{inj})^{min}(E, F) \stackrel{1}{=} ({\cal A}^{min})^{inj}(E, F)$
\end{center}
\end{prop}
{\noindent}{\sc PROOF}: It is sufficient to prove the first isometric equality
(for ${\cal A}^{sur}$) since the second one can be proved analogous.
Let $\oid{B}{}{} \in \{\oid{A}{}{},
\oid{A}{}{\displaystyle\ast\scriptstyle\triangle dd}\}$.
By assumption \oid{B}{}{} is a normed operator ideal and therefore 
$\oid{B}{}{\displaystyle\ast\displaystyle\ast} \stackrel{1}{=}
{\cal B}^{max}$ (see \cite{p1}, 9.3.1). Using known hull operations (cf.\cite{p1},
8.7 and 9.3) it follows that
\begin{center}
$({\cal B}^{sur})^{min} \stackrel{1}{=} 
({\cal B}^{\displaystyle{\ast\ast}\scriptstyle sur})^{min}
\stackrel{1}{=} ({\cal A}^{sur})^{min}$
\end{center}
and
\begin{center}
$({\cal B}^{min})^{sur} \stackrel{1}{=} 
({\cal B}^{\displaystyle{\ast\ast}\scriptstyle min})^{sur}
\stackrel{1}{=} ({\cal A}^{min})^{sur}$
\end{center}
In particular these equalities are true for $\oid{B}{}{}:= 
\oid{A}{}{\displaystyle\ast\scriptstyle\triangle dd}$.
Since \oid{B}{}{} is accessible (in particular right-accessible), the claim
follows by (\cite{df}, 25.11).${}_{\displaystyle\Box}$ \\

{\noindent}So far we have seen that conjugates of maximal Banach ideals
play a key role in the investigation of accessibility. Their appropriateness
will be strenghtened by illuminating accessibility via
a calculus derived from specific quotient ideals which are canonical extensions
of conjugate ideals and which appear in a natural matter by theorem 3.1.\\
Let (\oid{A}{}{}, \idn{A}{}{}) be a $p$-Banach ideal ($0 < p \leq 1$).\\

We put:\footnote{
In \cite{oe1}, $\oid{I}{}{} \circ \oid{A}{}{-1}$ was abbreviated as 
\oid{A}{}{\displaystyle\varepsilon} and $\oid{A}{}{-1} \circ \oid{I}{}{}$ as 
\oid{A}{\displaystyle\varepsilon}{}.}
\begin{center}
$(\oid{A}{}{\dashv}, \idn{A}{}{\dashv}) := (\oid{I}{}{} \circ \oid{A}{}{-1}, 
\idn{I}{}{} \circ \idn{A}{}{-1})$
\hspace{0.15cm}and\hspace{0.15cm}
$(\oid{A}{}{\vdash}, \idn{A}{}{\vdash}) := (\oid{A}{}{-1} \circ 
\oid{I}{}{}, \idn{A}{}{-1} \circ \idn{I}{}{})$
\end{center}
and omit the proof of the simple but useful

\begin{lemma}
\begin{itemize}
\item[(i)] \idrel{A}{\subseteq}{A^{\vdash\dashv}} and 
           \idrel{A}{\subseteq}{A^{\dashv\vdash}}
\item[(ii)] $\oid{A^{\scriptstyle\triangle}}{}{} \stackrel{1}{\subseteq}
            \oid{A^{\vdash}}{}{} \stackrel{1}{\subseteq}
            \oid{A^{\displaystyle\ast}}{}{} $
            and $\oid{A^{\scriptstyle\triangle}}{}{} \stackrel{1}{\subseteq}
            \oid{A^{\dashv}}{}{} \stackrel{1}{\subseteq}
            \oid{A^{\displaystyle\ast}}{}{} $
\item[(iii)] if \idrel{A}{\subseteq}{B} then \idrel{B^{\vdash}}{\subseteq}{A^{\vdash}} 
             and \idrel{B^{\dashv}}{\subseteq}{A^{\dashv}}.          
\end{itemize}
\end{lemma}

{\noindent}For $p=1$, theorem 3.1 therefore implies that \oid{A}{}{} is right-accessible 
if and only if 
\idrel{A^{\displaystyle\ast}}{=}{A^{\dashv}}. If 
\oid{A}{}{\displaystyle\ast\scriptstyle\triangle}
is left-accessible or \oid{A}{}{} maximal, then the left-accessibility of 
\oid{A}{}{} is equivalent to the statement 
\idrel{A^{\displaystyle\ast}}{=}{A^{\vdash}}.
Note that \idrel{A^{\scriptstyle\triangle}}{=}{A^{\vdash}} if \oid{A}{}{} is
injective and \idrel{A^{\scriptstyle\triangle}}{=}{A^{\dashv}} if \oid{A}{}{}
is surjective (see \cite{jo}, 2.6). Since \idn{A}{}{\scriptstyle\triangle}
and \idn{A}{}{\displaystyle\ast} coincide on the space of all elementary operators
it follows that always 
\idrel{A^{\scriptstyle\triangle\displaystyle\ast}}{=}{A^{\displaystyle\ast\displaystyle\ast}},
hence
\idrel{A^{\vdash\displaystyle\ast}}{=}{A^{\displaystyle\ast\displaystyle\ast}}.
Therefore, if \oid{A}{}{} is a {\it{maximal}} Banach ideal, then lemma 3.1 implies that 
$\oid{A}{}{} \stackrel{1}{=} 
\idrel{A^{\vdash\displaystyle\ast}}{=}{A^{\vdash\dashv}}$ and we have obtained the:
\begin{cor}
Let $(\oid{A}{}{}, \idn{A}{}{})$ be a maximal Banach ideal. Then 
$(\oid{A^{\vdash}}{}{}, \idn{A^{\vdash}}{}{})$ is right-accessible. If 
\oid{A^{\scriptstyle\triangle}}{}{} is left-accessible, then 
$(\oid{A^{\dashv}}{}{}, \idn{A^{\dashv}}{}{})$ is left-accessible.
\end{cor}

{\noindent}Now we turn our attention to {\it{injective}} maximal Banach ideals; 
in particular we are interested in aspects concerning the left-accessibility 
of those ideals.

\section{Totally accessible injective operator ideals}
Let  (\oid{A}{}{}, \idn{A}{}{}) be a maximal Banach ideal and 
(\oid{A}{}{inj}, \idn{A}{}{inj})
the injective hull of (\oid{A}{}{}, \idn{A}{}{}). Let $\alpha \sim \oid{A}{}{}$
be associated. Then $\alpha \setminus \sim \oid{A}{}{inj}$ and 
$\setminus\alpha^{\displaystyle\ast}  \sim \oid{A}{}{inj \displaystyle\ast}$.
Since $\setminus\alpha^{\displaystyle\ast} \sim  
\setminus\oid{A}{}{\displaystyle\ast}$
and $\oid{\setminus L}{}{} \stackrel{1}{=} \oid{L}{\infty}{}$ it follows that
\begin{center}
$\oid{A}{}{inj \displaystyle\ast} \stackrel{1}{=} \setminus\oid{A}{}{\displaystyle\ast}
\stackrel{1}{=} (\oid{A}{}{\displaystyle\ast} \circ \oid{\setminus L}{}{})^{reg}
\stackrel{1}{=} (\oid{A^{\displaystyle\ast}}{}{} \circ \oid{L}{\infty}{})^{reg} 
\sim \setminus\alpha^{\displaystyle\ast}$ 
\end{center}
is the adjoint of \oid{A}{}{inj}, hence a {\it{maximal Banach ideal}}
(see \cite{df}, 25.9). In particular we obtain $\oid{L}{\infty}{} \stackrel{1}{=}
(\oid{L}{\infty}{} \circ \oid{L}{\infty}{})^{reg}$ (since $\oid{P}{1}{} \stackrel{1}{=} 
\oid{L}{\infty}{\displaystyle\ast}$ is injective) and  $\oid{I}{}{} \stackrel{1}{=}
(\oid{P}{1}{} \circ \oid{L}{\infty}{})^{reg}$ (since $\oid{L}{\infty}{inj} \stackrel{1}{=} 
\oid{L}{}{}$ (cf. \cite{df}, 20.14)) which implies that both 
$(\oid{L}{\infty}{} \circ \oid{L}{\infty}{})^{reg}$ and 
$(\oid{P}{1}{} \circ \oid{L}{\infty}{})^{reg}$ are normed operator ideals - a fact
which is not obvious.
\begin{lemma}
Let $(\oid{A}{}{}, \idn{A}{}{})$ be a $p$-Banach ideal $(0 < p \leq 1)$ and
$(\oid{B}{}{}, \idn{B}{}{})$ be a $q$-Banach ideal $(0 < q \leq 1)$.
If $(\oid{A}{}{}, \idn{A}{}{}) \stackrel{1}{\subseteq} (\oid{A}{}{dd}, \idn{A}{}{dd})$
then
\begin{center}
$\oid{A}{}{} \circ \oid{B}{}{reg} \stackrel{1}{\subseteq} 
(\oid{A}{}{} \circ \oid{B}{}{})^{reg}$.

\end{center}
\end{lemma}

{\noindent}{\sc PROOF}: Let $E, F$ be Banach spaces, $\varepsilon > 0$ and $ T \in 
\oid{A}{}{} \circ \oid{B}{}{reg}(E, F)$. Then there are a Banach space $G$,
operators $R \in \oid{A}{}{}(G, F)$ and $S \in \oid{B}{}{reg}(E, G)$ such that
$T = RS$ and \\$\idn{A}{}{}(R) \idn{B}{}{reg}(S) < (1 + \varepsilon) (\idn{A}{}{} 
\circ \idn{B}{}{})^{reg}(T)$.
Hence $j_FT = R''j_GS \in \oid{A \circ B}{}{}(E, F'')$ and 
$\idn{A \circ B}{}{}(j_FT) \leq \idn{A}{}{}(R'') \idn{B}{}{reg}(S) <
(1 + \varepsilon) (\idn{A}{}{} \circ \idn{B}{}{})^{reg}(T)$. ${}_{\displaystyle\Box}$\\

{\noindent}Now we have prepared all tools to prove 
\begin{prop}
Let $(\oid{A}{}{}, \idn{A}{}{})$ be a maximal Banach ideal. Then the following
statements are equivalent:
\begin{itemize}
\item[(a)]  $\oid{A}{}{} \circ \oid{A^{\displaystyle\ast}}{}{} 
\stackrel{1}{\subseteq}
            \oid{P}{1}{}$
\item[(b)]  \oid{A}{}{inj} is totally accessible.
\end{itemize}
\end{prop}
{\noindent}{\sc PROOF}: Let (a) be valid. To prove (b), it is enough to show
that \oid{A}{}{inj \displaystyle\ast} is right accessible. Since \oid{P}{1}{}
is right-accessible, theorem 3.1 implies 
$\oid{A}{}{} \circ \oid{A^{\displaystyle\ast}}{}{} \circ \oid{L}{\infty}{} 
\stackrel{1}{\subseteq}
\oid{I}{}{}$ and hence 
$\oid{A^{\displaystyle\ast}}{}{} \circ \oid{L}{\infty}{} \stackrel{1}{\subseteq}
\oid{A}{}{\vdash}$. Since 
$\oid{L}{\infty}{} \stackrel{1}{=} 
(\oid{L}{\infty}{} \circ \oid{L}{\infty}{})^{reg}$, it follows by lemma 4.1
that $\oid{A^{\displaystyle\ast}}{}{} \circ \oid{L}{\infty}{} \stackrel{1}{\subseteq}
(\oid{A^{\displaystyle\ast}}{}{} \circ \oid{L}{\infty}{} \circ \oid{L}{\infty}{})^{reg} 
\stackrel{1}{\subseteq} (\oid{A}{}{\vdash} \circ \oid{L}{\infty}{})^{reg}$.
Hence
\begin{center}
$\oid{A}{}{inj \displaystyle\ast} \stackrel{1}{=} (\oid{A}{}{\vdash} 
\circ \oid{L}{\infty}{})^{reg}$.
\end{center}
Now we apply theorem 3.1 to the Banach ideal $(\oid{A}{}{\vdash} \circ 
\oid{L}{\infty}{})^{reg}$:
since \oid{A}{}{} is assumed to be a maximal Banach ideal and \oid{L}{\infty}{}
is right-accessible, $\oid{A}{}{\vdash} \circ \oid{L}{\infty}{} =: \oid{B}{}{}$
is also right-accessible. Hence $\oid{B^{\displaystyle\ast}}{}{} \circ \oid{B}{}{} 
\stackrel{1}{\subseteq}
\oid{I}{}{}$ and therefore 
\begin{center}
$(\oid{B}{}{reg})^{\displaystyle\ast} \circ \oid{B}{}{reg} \stackrel{1}{\subseteq}
\oid{B^{\displaystyle\ast}}{}{} \circ \oid{B}{}{reg} \stackrel{1}{\subseteq}
(\oid{B^{\displaystyle\ast}}{}{} \circ \oid{B}{}{})^{reg} 
\stackrel{1}{\subseteq} \oid{I}{}{}$.
\end{center}
In other words: $(\oid{B}{}{reg})^{\displaystyle\ast} \stackrel{1}{=}
(\oid{B}{}{reg})^{\dashv}$. Since $\oid{B}{}{reg} \stackrel{1}{=}
\oid{A}{}{inj \displaystyle\ast}$ is normed, theorem 3.1 implies that 
\oid{A}{}{inj \displaystyle\ast} is right-accesible.\\
Now let \oid{A}{}{inj} be totally accessible, hence left-accessible. Then
$\oid{A}{}{inj} \circ \oid{A}{}{inj \displaystyle\ast} \stackrel{1}{\subseteq}
\oid{I}{}{}$ and therefore 
$\oid{A}{}{} \circ \oid{A^{\displaystyle\ast}}{}{} \circ \oid{L}{\infty}{}
\stackrel{1}{\subseteq} 
\oid{A}{}{} \circ \oid{A}{}{inj \displaystyle\ast} \stackrel{1}{\subseteq}
\oid{I}{}{}$. Hence $\oid{A}{}{} \circ \oid{A^{\displaystyle\ast}}{}{} 
\stackrel{1}{\subseteq} \oid{L}{\infty}{\dashv} \stackrel{1}{=} 
\oid{L}{\infty}{\displaystyle\ast}
\stackrel{1}{=} \oid{P}{1}{}$. ${}_{\displaystyle\Box}$\\

{\noindent}We don't know if there exists a maximal Banach ideal \oid{A}{}{} such that 
\oid{A}{}{inj} is totally accessible and 
$\oid{A}{}{} \circ \oid{A}{}{\displaystyle\ast} \stackrel{1}{\not\subseteq}
\oid{I}{}{}$ (hence \oid{A}{}{} not left-accessible).
\begin{cor}
Let $(\oid{A}{}{}, \idn{A}{}{})$ be a maximal Banach ideal. Then 
$\oid{A}{}{} \circ \oid{L}{\infty}{}$ is left-accessible and 
\begin{center}
$\oid{A}{}{} \circ \oid{L}{\infty}{} \circ \oid{A^{\displaystyle\ast}}{}{}
\stackrel{1}{\subseteq} \oid{I}{}{}$.
\end{center}
If $\oid{A}{}{} \circ \oid{L}{\infty}{}$
is right-accessible, then \oid{A}{}{\displaystyle\ast inj} is totally accessible.
\end{cor}

{\noindent}{\sc PROOF}: Since $\oid{A}{}{} \circ \oid{L}{\infty}{} \stackrel{1}{\subseteq}
(\oid{A}{}{} \circ \oid{L}{\infty}{})^{reg} \stackrel{1}{=} \oid{\setminus A}{}{}$,
and $\idn{A}{}{} \circ \idn{L}{\infty}{}$ coincides with 
$(\idn{A}{}{} \circ \idn{L}{\infty}{})^{reg}$ on the space of all elementary operators,
$\oid{A}{}{} \circ \oid{L}{\infty}{}$ is left-accessible and 
\begin{center}
$\oid{A}{}{} \circ \oid{L}{\infty}{} \circ \oid{A^{\displaystyle\ast}}{}{} 
\stackrel{1}{\subseteq} (\oid{A}{}{} \circ \oid{L}{\infty}{})
\circ (\oid{A}{}{} \circ \oid{L}{\infty}{})^{\displaystyle\ast} \stackrel{1}{\subseteq} 
\oid{I}{}{}$.
\end{center}
Now let $\oid{A}{}{} \circ \oid{L}{\infty}{}$ be right-accessible. We have to show that
\begin{center}
$\oid{A^{\displaystyle\ast}}{}{} \circ \oid{A}{}{} \stackrel{1}{\subseteq}
\oid{P}{1}{}$.
\end{center}
But this follows from $\oid{A^{\displaystyle\ast}}{}{} \circ \oid{A}{}{} 
\circ \oid{L}{\infty}{} \stackrel{1}{\subseteq} 
(\oid{A}{}{} \circ \oid{L}{\infty}{})^{\displaystyle\ast} \circ
(\oid{A}{}{} \circ \oid{L}{\infty}{}) \stackrel{1}{\subseteq} \oid{I}{}{}$.
${}_{\displaystyle\Box}$\\

{\noindent}Now we will recognize that prop. 4.1 leads to interesting consequences concerning
the characterization of a class of injective maximal Banach ideals which
are totally accessible.\\ Since \oid{P}{1}{} is included, {\it{Grothendieck's inequality}}\/ in
operator form implies a non-trivial relation to \oid{L}{2}{} and \oid{L}{1}{}
in the following sense (with Grothendieck constant $K_G$):
\begin{lemma}
$\oid{L}{2}{} \circ \oid{L}{1}{} \subseteq \oid{P}{1}{}$ and
$\idn{P}{1}{}(T) \leq K_G \cdot (\idn{L}{2}{} \circ \idn{L}{1}{})(T)$ for all $T \in 
\oid{L}{2}{} \circ \oid{L}{1}{}$.
\end{lemma}
{\noindent}{\sc PROOF}: Let $E, F$ be Banach spaces, $\varepsilon > 0$ and $T \in 
\oid{L}{2}{} \circ \oid{L}{1}{}(E, F)$. Then there exists a Banach space G, 
operators $R \in \oid{L}{2}{}(G, F), S \in \oid{L}{1}{}(E, G)$ such that 
$T = RS$ and $\idn{L}{2}{}(R)\idn{L}{1}{}(S) < (\idn{L}{2}{} \circ \idn{L}{1}{})(T)$.
Since $\oid{L}{1}{} \stackrel{1}{=} \oid{L/}{}{}$ and $S \in \oid{L}{1}{}(E, G)$,
there exists a measure $\mu$, operators $W \in\oid{L}{}{}(L_1(\mu), G'')$ and
$Z \in\oid{L}{}{}(E, L_1(\mu))$ such that $j_GS = WZ$ and 
$\| W \| \| Z \| < (1 + \varepsilon) \idn{L}{1}{}(S)$.
Since $R'' \in \oid{L}{2}{}(G'', F'')$, there exists a Hilbert space H, operators
$U \in \oid{L}{}{}(H, F'')$ and $V \in \oid{L}{}{}(G'', H)$ such that $R'' = UV$ and
$\| U \|\| V \| < (1 + \varepsilon) \                                  
\idn{L}{2}{}(R)$. Hence $j_FT = R''j_GS = U(VW)Z$ and $VW \in 
\oid{L}{}{}(L_1(\mu), H)$.
Since $L_1(\mu)$ is \oid{L}{1,1}{g} - space, Grothendieck's inequality implies 
that $VW \in \oid{P}{1}{}(L_1(\mu), H)$ and $\idn{P}{1}{}(VW) \leq K_G \cdot
\| VW \|$ (cf. \cite{df}, 23.10).\\ 
Hence $j_FT \in \oid{P}{1}{}(E, F'')$
and $\idn{P}{1}{}(j_FT) \leq \| U \| \idn{P}{1}{}(VW)
\| Z \| \leq (1 + \varepsilon)^2\ K_G\ \idn{L}{2}{}(R)
\idn{L}{1}{}(S)\\ <(1 + \varepsilon)^3\ K_G\ (\idn{L}{2}{} \circ \idn{L}{1}{})(T)$.
Since \oid{P}{1}{} is regular, the claim follows.
${}_{\displaystyle\Box}$ \\

{\noindent}Now let (\oid{A}{}{}, \idn{A}{}{}) be a maximal Banach ideal such that 
$\oid{D}{2}{} \subseteq \oid{A}{}{} \subseteq \oid{L}{1}{}$ (since 
$\oid{P}{1}{}  \stackrel{1}{\subseteq} \oid{L}{2}{}  \stackrel{1}{=} \oid{L}{2}{d}$, 
it follows
that $\oid{D}{2}{} \stackrel{1}{=} \oid{L}{2}{\displaystyle\ast} 
\stackrel{1}{\subseteq} \oid{P}{1}{d \displaystyle\ast} \stackrel{1}{=}
\oid{L}{1}{}$. Hence the class of such ideals \oid{A}{}{} is not empty; consider e.g. 
\oid{P}{2}{d}). Then
$\oid{A^{\displaystyle\ast}}{}{} \subseteq \oid{L}{2}{}$ and therefore 
$\oid{A^{\displaystyle\ast}}{}{} \circ \oid{A}{}{} \subseteq 
\oid{L}{2}{} \circ \oid{L}{1}{} \subseteq \oid{P}{1}{}$ by lemma 4.2.\\
If $\oid{A^{\displaystyle\ast}}{}{} \circ \oid{A}{}{} \stackrel{1}{\subseteq}
\oid{P}{1}{}$, prop. 4.1 would imply that \oid{A}{}{\displaystyle\ast inj}
is totally accessible. In general we don't know if this is the case. However
there exists a beautiful "trick" to arrange $\idn{P}{1}{}(T) \leq 1 \cdot 
(\idn{A^{\displaystyle\ast}}{}{} \circ \idn{A}{}{})(T)$ for all $T \in
\oid{A^{\displaystyle\ast}}{}{} \circ \oid{A}{}{}$, which is given by
{\it{tensor stability}}.\\
Let $\gamma$ be a fixed tensor norm. Remember that a given quasi Banach ideal 
(\oid{A}{}{}, \idn{A}{}{}) is called {\it{$\gamma$-tensorstable}} 
(cf. \cite{cdr}, \cite{df}), if
\begin{center}
$S \tilde{\otimes}_\gamma T \in \oid{A}{}{}(E \tilde{\otimes}_\gamma F, 
G \tilde{\otimes}_\gamma H)$
for all $S \in \oid{A}{}{}(E, G), T \in \oid{A}{}{}(F, H)$.
\end{center}
In this case there is a constant $c \geq 1$ satisfying
\begin{center}
$\idn{A}{}{}(S \tilde{\otimes}_\gamma T) \leq c \idn{A}{}{}(S) \idn{A}{}{}(T)$.
\end{center}
If $c=1$, \oid{A}{}{} is called {\it{metrically $\gamma$-tensorstable}}.
If $c=1$ and the above inequality is an equality, then
\oid{A}{}{} is called {\it{strongly $\gamma$-tensorstable}}.\\
With the help of tensor stability, we will show how it is possible 
to improve the inequality 
$\oid{A^{\displaystyle\ast}}{}{} \circ \oid{A}{}{} \subseteq
\oid{P}{1}{}$\ to obtain
$\oid{A^{\displaystyle\ast}}{}{} \circ \oid{A}{}{} \stackrel{1}{\subseteq}
\oid{P}{1}{}$\ and turn our attention to  
\begin{theorem}
Let $(\oid{A}{}{}, \idn{A}{}{})$ be a maximal Banach ideal. If
\begin{itemize}
\item[(i)]   $\oid{D}{2}{} \subseteq \oid{A}{}{} \subseteq \oid{L}{1}{}$ and 
\item[(ii)]  both $(\oid{A^{\displaystyle\ast}}{}{}, \idn{A^{\displaystyle\ast}}{}{})$
and $(\oid{A}{}{}, \idn{A}{}{})$ are metrically $\varepsilon$-tensorstable
\end{itemize}
then $(\oid{A}{}{\displaystyle\ast inj}, \idn{A}{}{\displaystyle\ast inj})$
is totally accessible.
\end{theorem}

{\noindent}{\sc PROOF}: Since  $\oid{D}{2}{} \subseteq \oid{A}{}{} \subseteq 
\oid{L}{1}{}$,\
the adjoint \oid{A^{\displaystyle\ast}}{}{} is contained in \oid{L}{2}{}. Hence
there exist constants $c \geq 0$ and $c^{\ast} \geq 0$ such that 
$\idn{L}{2}{}(R) \leq c^{\ast} \idn{A^{\displaystyle\ast}}{}{}(R)$ for all 
$R \in \oid{A}{}{\displaystyle\ast}$ and $\idn{L}{1}{}(S) \leq c \idn{A}{}{}(S)$
for all $S \in \oid{A}{}{}$. Let $E, F$ be Banach spaces, $\varepsilon >0$
and $T \in \oid{A}{}{\displaystyle\ast} 
\circ \oid{A}{}{}(E, F)$. We must show that $T \in \oid{P}{1}{}(E, F)$ and
$\idn{P}{1}{}(T) \leq 1 \cdot (\idn{A}{}{\displaystyle\ast} \circ \idn{A}{}{})(T)$.\\
By the previous considerations, there exists a Banach space $D$, operators
$R \in \oid{A}{}{\displaystyle\ast}(D, F)$ and $S \in \oid{A}{}{}(E, D)$
such that $T = RS \in \oid{L}{2}{} \circ \oid{L}{1}{}$ 
and $\idn{L}{2}{}(R)\ \idn{L}{1}{}(S) < cc^\ast 
(1 + \varepsilon) (\idn{A}{}{\displaystyle\ast} \circ \idn{A}{}{})(T)$.
Lemma 4.2 now implies that $T = RS \in \oid{P}{1}{}(E, F)$ and 
\begin{center}
$\idn{P}{1}{}(T) \leq K_G \idn{L}{2}{}(R) \idn{L}{1}{}(S) <
(1 + \varepsilon) K_G cc^\ast (\idn{A}{}{\displaystyle\ast} \circ \idn{A}{}{})(T)$.
\end{center}
At this point the improvement of this norm estimation will be realized by the assumed
metric $\varepsilon$-tensor stability of \oid{A^{\displaystyle\ast}}{}{}
and \oid{A}{}{} which implies in particular that $\oid{A}{}{\displaystyle\ast} \circ 
\oid{A}{}{}$ is
is metrically $\varepsilon$-tensorstable (cf. \cite{df}, 34.4). Since
\oid{P}{1}{} even is {\it{strongly}} $\varepsilon$-tensorstable
(see \cite{df}, 34.5), it follows that
\begin{eqnarray}
\idn{P}{1}{}{}(T)^2 & = & \idn{P}{1}{}(T \tilde{\otimes}_\varepsilon T)
\nonumber \\
& \leq & (1 + \varepsilon) K_G cc^\ast  
(\idn{A}{}{\displaystyle\ast} \circ \idn{A}{}{})(T \tilde{\otimes}_\varepsilon T)
\nonumber \\
& \leq & (1 + \varepsilon) K_G cc^\ast  
(\idn{A}{}{\displaystyle\ast} \circ \idn{A}{}{})(T)^2. \nonumber
\end{eqnarray}
Hence:\ $\idn{P}{1}{}{}(T) \leq ((1 + \varepsilon) K_G  cc^\ast)^{1/2}\
(\idn{A}{}{\displaystyle\ast} \circ \idn{A}{}{})(T)$, 
and an obvious induction argument implies:
\begin{center}
$\forall n \in {\bf N}:\ \idn{P}{1}{}{}(T) \leq ((1 + \varepsilon) K_G cc^\ast)^{1/2^n}\
(\idn{A}{}{\displaystyle\ast} \circ \idn{A}{}{})(T)$.
\end{center}
$n \rightarrow \infty$ now yields the desired improved norm estimation, and
the proof is finished.
${}_{\displaystyle\Box}$ \\

{\noindent}Next we will show that the statement of theorem 4.1 remains valid for 
arbitrary finitely generated tensor norms (not only for the injective
tensor norm $\varepsilon$) if we assume a (slight) restriction of the tensor
stability condition - with a completely different proof than the previous one.\\
So, let (\oid{A}{}{}, \idn{A}{}{}) be a maximal Banach ideal with 
$\oid{D}{2}{} \subseteq \oid{A}{}{} \subseteq \oid{L}{1}{}$. Then $\oid{\setminus A}{}{}
\subseteq \oid{\setminus L}{1}{} \stackrel{1}{=} 
(\oid{P}{1}{sur})^{\displaystyle\ast}$. By Grothendieck's inequality, 
$\oid{L}{2}{} \subseteq \oid{P}{1}{sur}$ and $\idn{P}{1}{sur}(T) \leq 
K_G \idn{L}{2}{}(T)$ for all $T \in \oid{L}{2}{}$ (cf. \cite{df}, 20.17)
Hence $\oid{\setminus A}{}{} \subseteq \oid{D}{2}{}$ and there is a constant
$c \geq 0$  such that $\idn{D}{2}{}(S) \leq  cK_G\idn{\setminus A}{}{}(S)$ for
all $S \in \oid{\setminus A}{}{}$. Since \oid{L}{2}{} is injective, it follows
that $(\oid{\setminus A}{}{})^{\displaystyle\ast} \subseteq \oid{L}{2}{}$, and there exists
a constant $c^{\ast}$ such that $\idn{L}{2}{}(R) \leq c^{\ast}
(\idn{\setminus A}{}{})^{\displaystyle\ast}(R)$  for all $R \in 
(\oid{\setminus A}{}{})^{\displaystyle\ast}$.\\
Since \oid{D}{2}{} is right-accessible, hence 
$\oid{L}{2}{} \circ \oid{D}{2}{} \stackrel{1}{\subseteq} \oid{I}{}{}$,
we have obtained the following statement which is of own interest:
\begin{lemma}
Let $(\oid{A}{}{}, \idn{A}{}{})$ be a maximal Banach ideal such that 
$\oid{D}{2}{} \subseteq \oid{A}{}{} \subseteq \oid{L}{1}{}$. Then
\begin{center}
$(\oid{\setminus A}{}{})^{\displaystyle\ast} \circ \oid{\setminus A}{}{}
\subseteq \oid{I}{}{}$ 
\end{center}
and there exist constants $c, c^{\ast}$ such that
\begin{center}
$\idn{I}{}{}(T) \leq cc^{\ast}K_G \hspace{0.1cm} 
((\idn{\setminus A}{}{})^{\displaystyle\ast} \circ \idn{\setminus A}{}{})(T)$ 
\end{center}
for all $T \in (\oid{\setminus A}{}{})^{\displaystyle\ast} \circ \oid{\setminus A}{}{}.
\hspace{0.5cm}
{}_{\displaystyle\Box}$
\end{lemma} 

{\noindent}Now let in addition $\gamma$ be an arbitrary finitely generated 
tensor norm and assume that $(\oid{\setminus A}{}{})^{\displaystyle\ast}$ as well as 
\oid{\setminus A}{}{} are metrically $\gamma$-tensorstable.\\
Let $(M, F) \in$ FIN $\times$ BAN and $U \in  
(\oid{\setminus A}{}{})^{\displaystyle\ast \scriptstyle\vdash}(M, F)$.
Then $U$ and $U \tilde{\otimes}_\gamma U$ are finite operators
and lemma 4.3. implies that 
\begin{center}
$(\idn{\setminus A}{}{})^{\displaystyle\ast \scriptstyle\vdash}(U) \leq
cc^{\ast}K_G\idn{\setminus A}{}{}(U)$
\end{center}
This estimation now can be improved as follows:\\
Let $\varepsilon > 0$. Then there is a Banach space $D$, an operator $V \in
(\oid{\setminus A}{}{})^{\displaystyle\ast}(F, D)$ with 
$(\idn{\setminus A}{}{})^{\displaystyle\ast}(V) = 1$ such that
$(\idn{\setminus A}{}{})^{\displaystyle\ast \scriptstyle\vdash}(U) < 
(1 + \varepsilon)\idn{I}{}{}(VU)$. Since $\oid{L}{}{} \stackrel{1}{=} 
\oid{I}{}{\displaystyle\ast}$ as well as \oid{I}{}{} are metrically $\gamma$-
tensorstable (see \cite{df}, 34.5), \oid{I}{}{} even is strongly $\gamma$-
tensorstable (see \cite{df}, 34.2). Hence  
\begin{eqnarray}
(\idn{\setminus A}{}{})^{\displaystyle\ast \scriptstyle\vdash}(U)^2 & < &
(1 + \varepsilon)^2\idn{I}{}{}(VU)^2 \nonumber\\
& \leq & (1 + \varepsilon)^2\idn{I}{}{}((V \tilde{\otimes}_\gamma V) \circ
(U \tilde{\otimes}_\gamma U)) \nonumber\\
& \leq & (1 + \varepsilon)^2 
(\idn{\setminus A}{}{})^{\displaystyle\ast}(V \tilde{\otimes}_\gamma V)\
(\idn{\setminus A}{}{})^{\displaystyle\ast \scriptstyle\vdash}(U \tilde{\otimes}_
\gamma U) 
\nonumber \\
& \leq & (1 + \varepsilon)^2 
(\idn{\setminus A}{}{})^{\displaystyle\ast \scriptstyle\vdash}(U \tilde{\otimes}_
\gamma U) 
\nonumber 
\end{eqnarray}
(The last inequality follows by the metric $\gamma$-tensor stability
of $(\oid{\setminus A}{}{})^{\displaystyle\ast})$.
Since \oid{\setminus A}{}{} also is metrically $\gamma$-tensorstable, we obtain
\begin{center}
$(\idn{\setminus A}{}{})^{\displaystyle\ast \scriptstyle\vdash}(U)^2 \leq
(\idn{\setminus A}{}{})^{\displaystyle\ast \scriptstyle\vdash}(U \tilde{\otimes}_
\gamma U) 
\leq
cc^{\ast}K_G\ \idn{\setminus A}{}{}(U \tilde{\otimes}_\gamma U) \leq
cc^{\ast}K_G\ \idn{\setminus A}{}{}(U)^2$.
\end{center}
Hence, induction implies 
$(\idn{\setminus A}{}{})^{\displaystyle\ast \scriptstyle\vdash}(U) \leq
\idn{\setminus A}{}{}(U)$, and since 
$(\oid{\setminus A}{}{})^{\displaystyle\ast \scriptstyle\vdash}$ is 
right-accessible, we have proved:

\begin{theorem}
Let $(\oid{A}{}{}, \idn{A}{}{})$ be a maximal Banach ideal and $\gamma$ be a
finitely generated tensor norm. If
\begin{itemize}
\item[(i)]   $\oid{D}{2}{} \subseteq \oid{A}{}{} \subseteq \oid{L}{1}{}$ and 
\item[(ii)]  both $((\oid{\setminus A}{}{})^{\displaystyle\ast}, 
(\idn{\setminus A}{}{})^{\displaystyle\ast})$
and $(\oid{\setminus A}{}{}, \idn{\setminus A}{}{})$ are metrically $\gamma$-tensorstable
\end{itemize}
then $(\oid{A}{}{\displaystyle\ast inj}, \idn{A}{}{\displaystyle\ast inj})
\stackrel{1}{=} ((\oid{\setminus A}{}{})^{\displaystyle\ast}, 
(\idn{\setminus A}{}{})^{\displaystyle\ast})$
is totally accessible.
\hspace{0.5cm} ${}_{\displaystyle\Box}$
\end{theorem}

{\noindent}To finish this paper we will give now interesting examples which also 
show that {\it{normed products}} of maximal Banach ideals will be of
crucial significance concerning the investigation of accessibility. In particular
we give a partial answer to a question of A. Defant and K. Floret whether
the ideal (\oid{L}{\infty}{}, \idn{L}{\infty}{}) is totally accessible or
not (see \cite{df}, 21.12).

\section{On normed products of operator ideals and accessibility}
Let (\oid{A}{}{}, \idn{A}{}{}) be a maximal Banach ideal such that
$(\oid{A}{}{}, \idn{A}{}{}) \subseteq (\oid{L}{2}{}, \idn{L}{2}{})$.\\
Then $\oid{D}{2}{} \subseteq \oid{A^{\displaystyle\ast}}{}{}$, and related to 
the foregoing results, it is a natural question to ask for further properties of \oid{A}{}{},
which even imply that $\oid{D}{2}{}  \subseteq \oid{A}{}{\displaystyle\ast} \subseteq
\oid{L}{1}{}$. In order to arrange this inclusion, we consider now
the product ideal $\oid{L}{2}{} \circ \oid{A}{}{}$:
\begin{lemma}
Let $(\oid{A}{}{}, \idn{A}{}{})$ be a $p$-Banach ideal $(0 < p \leq 1)$.
Then $(\oid{L}{2}{} \circ \oid{A}{}{}, \idn{L}{2}{} \circ \idn{A}{}{})$ always is
an injective $\frac{p}{1+p}$-Banach ideal. In particular it is 
right-accessible and regular.
\end{lemma}

{\noindent}{\sc PROOF}: Let $T \in (\oid{L}{2}{} \circ \oid{A}{}{})^{inj}(E, F)$ and
$\varepsilon > 0$. Then there are a Banach space $G$, operators $R \in 
\oid{L}{2}{}(G, F^{\infty})$
and $S\in \oid{A}{}{}(E, G)$ such that $J_FT = RS$ and $\idn{L}{2}{}(R) \idn{A}{}{}(S) < 
(1 + \varepsilon) \idn{L}{2}{} \circ \idn{A}{}{}(J_FT)$. Let $H$ be a 
Hilbert space, $V \in \oid{L}{}{}(H, F^{\infty})$ and $W \in \oid{L}{}{}(G, H)$
such that $R = VW$ and $\| V \| \| W \| < (1 + \varepsilon) 
\idn{L}{2}{}(R)$. Let $C$ be the (closed) range of $J_F : F \stackrel{1}{\hookrightarrow} 
F^{\infty}$. 
Then $H_0 : = V^{-1}(C)$ is a closed subspace of $H$, 
and consequently there exists a projection $P \in \oid{L}{}{}(H, H)$
from $H$ onto $H_0$ such that the closure of the range of $VP$ is contained 
in $C$. Since the range of $WS$ is contained in $H_0$, it follows that 
$WS = PWS$.\\
Now let $\gamma : C \longrightarrow F$ be defined canonically and let
$\gamma_0$ be the restriction of $\gamma$ to $C_0$, where $C_0$ is
the closure of $VP(H)$. Let $B : = \gamma_0 Z$,
with $Z : H \longrightarrow C_0, z \mapsto VPz$, and let $D : = WS$.
Then, $B \in \oid{L}{2}{}(H, F)$ and $\idn{L}{2}{}(B) \leq 
\| V \|$, $D \in \oid{A}{}{}(E, H)$ and  $\idn{A}{}{}(D) \leq 
\| W \| \idn{A}{}{}(S)$. In accordance with the construction,
\begin{center}
$BDx = \gamma(VPWSx) = \gamma(VWSx) = \gamma(RSx) = Tx$\hspace{0.3cm} for all 
$x \in E$. 
\end{center}
Hence $T = BD \in 
\oid{L}{2}{} \circ \oid{A}{}{}(E, F)$ and 
$\idn{L}{2}{}(B) \idn{A}{}{}(D) < (1 + \varepsilon)^2 
(\idn{L}{2}{} \circ \idn{A}{}{})^{inj}(T)$.
Injective $p$-Banach ideals are always right-accessible, and since $F^{\infty}$
has the metric extension property, they are also regular.
${}_{\displaystyle\Box}$\\

{\noindent}Which ideals \oid{A}{}{} imply now the non-normability of
$\oid{L}{2}{} \circ \oid{A}{}{}$? One answer is given by 

\begin{prop}
Let $(\oid{A}{}{}, \idn{A}{}{})$ be a maximal Banach ideal such that
$(\oid{L}{2}{} \circ \oid{A}{}{}, \idn{L}{2}{} \circ \idn{A}{}{})$
is normed. Then $(\oid{A}{}{\displaystyle\ast}, \idn{A}{}{\displaystyle\ast}) 
\subseteq (\oid{L}{\infty}{}, \idn{L}{\infty}{})$.
\end{prop}

{\noindent}{\sc PROOF}: Since $\oid{L}{2}{} \circ \oid{A}{}{}$ is normed and 
injective (in particular regular),
and as a product of two ultrastable Banach ideals again ultrastable
(cf. \cite{d}, 3.4.5), it follows that 
\begin{center}
$\oid{N}{}{} \stackrel{1}{\subseteq} \oid{L}{2}{} \circ \oid{A}{}{} 
\stackrel{1}{=} (\oid{L}{2}{} \circ \oid{A}{}{})^{reg} \stackrel{1}{=}
(\oid{L}{2}{} \circ \oid{A}{}{})^{max}$. 
\end{center}
The last isometric identity is implied by 
(\cite{p1}, 8.8.6). Hence $\oid{I}{}{} \stackrel{1}{=} \oid{N}{}{max} 
\stackrel{1}{\subseteq} \oid{L}{2}{} \circ \oid{A}{}{}$, and the injectivity further implies
that $\oid{P}{1}{} \stackrel{1}{\subseteq} 
\oid{L}{2}{} \circ \oid{A}{}{}$ and therefore
\begin{center}
$\oid{A}{}{\displaystyle\ast} \stackrel{1}{\subseteq}
(\oid{L}{2}{} \circ \oid{A}{}{})^{\displaystyle\ast} 
\stackrel{1}{\subseteq} \oid{L}{\infty}{}$.
${}_{\displaystyle\Box}$\\
\end{center}

{\noindent}Combining the preceeding considerations with theorem 4.1 and 4.2,
lead to the somehow surprising 

\begin{cor}
Let $(\oid{A}{}{}, \idn{A}{}{})$ be a maximal Banach ideal such that
$(\oid{A}{}{}, \idn{A}{}{}) \subseteq (\oid{L}{2}{}, \idn{L}{2}{})$.\\
Let $(\oid{A^{\displaystyle\ast}}{}{}, \idn{A^{\displaystyle\ast}}{}{})$ as
well as $(\oid{A}{}{}, \idn{A}{}{})$ be metrically $\varepsilon$-tensorstable, or
let $(\oid{\setminus A}{}{\displaystyle\ast}, \idn{\setminus A}{}{\displaystyle\ast})$
and $(\oid{A}{}{inj}, \idn{A}{}{inj})$ be metrically $\gamma$-tensorstable
with respect to a given finitely generated tensor norm $\gamma$.\\
If $(\oid{L}{2}{} \circ \oid{A}{}{d}, \idn{L}{2}{} \circ \idn{A}{}{d})$
is normed, then $(\oid{A}{}{inj}, \idn{A}{}{inj})$ is totally accessible.
\end{cor}

{\noindent}{\sc PROOF}: Let $\oid{B}{}{} : = \oid{A}{}{d}$. Then
\begin{center}
$\oid{D}{2}{} \subseteq \oid{A}{}{\displaystyle\ast}
\stackrel{1}{=} \oid{B}{}{\displaystyle\ast d} 
\stackrel{1}{\subseteq} \oid{L}{\infty}{\displaystyle d} \stackrel{1}{=}
\oid{L}{1}{}$.
\end{center}
Now, theorem 4.1 and 4.2 - applied to $\oid{A}{}{\displaystyle\ast}$ - yield the claim.
${}_{\displaystyle\Box}$\\

\begin{cor}
Let $(\oid{A}{}{}, \idn{A}{}{})$ be a maximal injective Banach ideal such that
$(\oid{A}{}{}, \idn{A}{}{}) \subseteq (\oid{L}{2}{}, \idn{L}{2}{})$.\\
Let $(\oid{A^{\displaystyle\ast}}{}{}, \idn{A^{\displaystyle\ast}}{}{})$ as
well as $(\oid{A}{}{}, \idn{A}{}{})$ be metrically $\gamma$-tensorstable
with respect to a given finitely generated tensor norm $\gamma$.
If $(\oid{L}{2}{} \circ \oid{A}{}{d}, \idn{L}{2}{} \circ \idn{A}{}{d})$ 
is normed, then $(\oid{A}{}{}, \idn{A}{}{})$ is totally accessible.
\end{cor}

\begin{prop}
Let $(\oid{A}{}{}, \idn{A}{}{})$ be a maximal Banach ideal such that
$(\oid{A^{\displaystyle\ast}}{}{}, \idn{A^{\displaystyle\ast}}{}{})$ is 
$\gamma$-tensorstable for some injective tensor norm $\gamma$.
If $(\oid{L}{2}{} \circ \oid{A}{}{d}, \idn{L}{2}{} \circ \idn{A}{}{d})$
is normed, than there is no infinite dimensional Banach space $E$ such that
$id_E \in \oid{A}{}{\displaystyle\ast}$.
\end{prop} 
{\noindent}{\sc PROOF}: By (\cite{df}, 23.3), the ideals \oid{L}{1}{},
\oid{L}{2}{}, and \oid{L}{\infty}{} are mutually uncomparable. Since 
$\gamma$ is an injective tensor norm, an infinite dimensional Banach space
in space$(\oid{A}{}{\displaystyle\ast})$ would lead to $\oid{L}{\infty}{} \subseteq 
\oid{A}{}{\displaystyle\ast}$ (cf. \cite{df}, 34.7). 
On the other hand, (since $\oid{L}{2}{} \circ \oid{A}{}{d}$ is normed)
the previous calculations show that
$\oid{A}{}{\displaystyle\ast} \stackrel{1}{\subseteq} \oid{L}{1}{}$.
Hence
$\oid{L}{\infty}{} \subseteq \oid{A}{}{\displaystyle\ast} 
\stackrel{1}{\subseteq}
\oid{L}{1}{}$, which is a contradiction.
${}_{\displaystyle\Box}$\\

{\noindent}{\bf{Remark:}} Let $\oid{A}{}{} \subseteq \oid{D}{2}{}$. Then
$\oid{A}{}{} \subseteq \oid{L}{2}{}$, and $\oid{L}{2}{} \circ \oid{A}{}{d}$
is a {\it{trace ideal}} (cf. \cite{d}, 4.4). In particular 
$\oid{L}{2}{} \circ \oid{A}{}{d}$ is not normed.\\

{\noindent}So far we have seen that there is an intimate relation between normed products
of operator ideals and accessibility-conditions. We will finish this paper
with another example which shows again, that normed products of operator ideals have an 
impact on accessibility.

\begin{prop}
Let $(\oid{A}{}{}, \idn{A}{}{})$ be an injective, maximal Banach ideal, which
is totally accessible. If 
$((\oid{A}{}{\displaystyle\ast} \circ \oid{A}{}{})^{reg}, 
(\idn{A}{}{\displaystyle\ast} \circ \idn{A}{}{})^{reg})$ 
is normed, then
$(\oid{A}{}{\displaystyle\ast}, \idn{A}{}{\displaystyle\ast})$ is not
totally accessible.
\end{prop}

{\noindent}{\sc PROOF}: By assumption, both \oid{A}{}{} and \oid{A}{}{\displaystyle\ast}
are maximal, hence ultrastable, so is their
product (see \cite{d}, 3.4.5), and it follows that
\begin{center} 
$(\oid{A}{}{\displaystyle\ast} \circ \oid{A}{}{})^{reg} \stackrel{1}{=}
(\oid{A}{}{\displaystyle\ast} \circ \oid{A}{}{})^{max}$.
\end{center}
Since $(\oid{A}{}{\displaystyle\ast} \circ \oid{A}{}{})^{reg}$ is assumed to be normed,
and \oid{A}{}{} is right-accessible, we obtain
\begin{center}
$\oid{I}{}{} \stackrel{1}{=} \oid{N}{}{max} \stackrel{1}{\subseteq}
(\oid{A}{}{\displaystyle\ast} \circ \oid{A}{}{})^{reg} \stackrel{1}{\subseteq}
\oid{I}{}{}$.
\end{center}
Since \oid{A}{}{} is injective  and \oid{A}{}{\displaystyle\ast} regular, an easy
calculation shows that $\oid{A}{}{\displaystyle\ast} \circ 
\oid{A}{}{}$ is regular as well, and therefore it follows that
\begin{center}
$\oid{A}{}{\displaystyle\ast} \circ \oid{A}{}{} \stackrel{1}{=}
\oid{I}{}{}$.
\end{center}
Now, assume that \oid{A}{}{\displaystyle\ast} is totally accessible.
Then, by the injectivity of the totally accessible ideal \oid{A}{}{},\hspace{0.1cm} 
$\oid{A}{}{\displaystyle\ast} \circ \oid{A}{}{}$ is also {\it{totally}} accessible
(see \cite{df}, 21.4), which is a contradiction, since \oid{I}{}{} is
not totally accessible.
${}_{\displaystyle\Box}$\\

\begin{cor}
If $(\oid{L}{\infty}{} \circ \oid{P}{1}{})^{reg}$ is normed, then \oid{L}{\infty}{}
is not totally accessible.
\end{cor}
Note, that  $(\oid{P}{1}{} \circ \oid{L}{\infty}{})^{reg} \stackrel{1}{=} 
\oid{I}{}{}$ is normed, as it was shown in section 4.                   

\section{Questions and open problems}
\begin{itemize}
\item  Is the conjugate of a maximal Banach ideal always left-accessible?
       (conjecture: no)  
\item  Let $(\oid{A}{}{}, \idn{A}{}{})$ be a maximal Banach ideal. Does this
       even imply the validity of the \oid{A}{}{}-local principle of 
       reflexivity?        
       (conjecture: no)
\item  What relations exist between tensorstable operator ideals, normed products of
       operator ideals and accessibility? How far, {\it{trace ideals}} are
       involved?               
\item  Assume, \oid{A}{}{} is a maximal Banach ideal such that 
       $\oid{A}{}{\displaystyle\ast} \stackrel{1}{\subseteq} \oid{L}{\infty}{}$.
       Does this condition even imply that 
       $\oid{L}{2}{} \circ \oid{A}{}{}$
       is a {\it{normed}} ideal?
       (The converse implication is true (see 5.1)).
\item  Is it possible to maintain the statement of corollary 5.2 {\it{without}}
       the property of $\oid{L}{2}{} \circ \oid{A}{}{d}$ being normed?
\item  In general it seems to be more easy (via trace ideals) to prove 
       that the product of two Banach ideals is {\it{not}} normed. 
       Find criteria, which show the normability.
\end{itemize}

{\noindent}Frank Oertel\\
Swiss Reinsurance Company\\
Dpt.: PM-PH (Development and Support)\\
Mythenquai 50/60\\
CH - 8022 Zurich\\
SWITZERLAND\\

{\noindent}E-mail: chvs95gq@ibmmail.com\\
Tel.: +41-1-285-3688\\

\end{document}